\documentclass[letterpaper, 10 pt, conference]{ieeeconf}  

\IEEEoverridecommandlockouts                              
\usepackage{graphicx}
\usepackage{epstopdf}
\usepackage{amsmath}
\usepackage{multicol}
\usepackage{stfloats}
\usepackage{amsfonts}
\usepackage{mathrsfs}
\usepackage{color}

\newtheorem{lemma}{Lemma}
\newtheorem{proposition}{Proposition}

\title{\LARGE \bf Distributed estimation and control for preserving  formation rigidity for mobile robot teams}

\author{Zhiyong Sun, Changbin Yu and Brian D. O. Anderson
\thanks{Zhiyong Sun is  with Shandong Computer Science Center (SCSC), Jinan, China;  Brian D. O. Anderson  was a visiting expert with SCSC.} \thanks{Zhiyong Sun, Changbin Yu and Brian D. O. Anderson are with National ICT Australia and
Research School of Engineering, The Australian National University, Canberra ACT
0200, Australia.   ({\tt\small \{zhiyong.sun, brad.yu, brian.anderson\}@anu.edu.au})}%
\thanks{}%
}

\begin{document}

\maketitle
\thispagestyle{empty}
\pagestyle{empty}


\begin{abstract}
Inspired by the concept of network algebraic connectivity, we adopt an extended notion named rigidity preservation index to characterize the rigidity property for a formation framework. A gradient based controller is proposed to ensure the rigidity preservation of multi-robot networks in an unknown environment, while the rigidity metric can be maximized over time during robots' motions. In order to implement the controller in a distributed manner, a distributed inverse power iteration algorithm is developed which allows each robot to estimate the global rigidity index information. Simulation results are provided to demonstrate the effectiveness of the estimation and control scheme.
\end{abstract}
\section{Introduction}
Networked mobile robots have the desirable capacity of performing spatially distributed tasks like  large area surveillance, underwater exploration, target detection, etc, while these tasks generally cannot be achieved by single robot. To deploy and coordinate a group of mobile robots in a complicated or even perhaps uncertain environment, a fundamental control problem  is how to ensure the connectivity of mobile robots under communication constraints. In fact, network connectivity maintenance and control is always a critical issue for cooperative robotics  and it has received considerable attention in recent years \cite{stump2008connectivity},  \cite{zavlanos2011graph}, \cite{sabattin2013}.

In robot coordination control design,  another closely related concept which has  equal importance to network connectivity is formation rigidity. It has been shown that rigid graph theory plays a key role in analyzing the formation performance \cite{anderson2008rigid} and network localization \cite{aspnes2006theory}. Formation rigidity can  be particularly important in the formation shape control for multi-robot systems, which enables the distance-based control design without a global coordinate or a centralized control. The  favorable property of a rigid framework, therefore,  motivates us to consider the problem of preserving formation rigidity when coordinating robot teams with planned motions. This can be regarded as an extension of the network connectivity control and is different from the formation shape control. In fact, to preserve only the formation rigidity also renders robot teams some level of flexibility for performing other tasks without a strict geometric shape requirement.

   The preservation of  the formation rigidity for mobile robot teams, rather than the mere maintenance of network connectivity, possesses  several advantages.     Compared with a connected robot team, a rigid dynamic formation generally has a high level of robustness on recovering from link or agent failures and it is also more convenient to perform a desired formation reconfiguration \cite{anderson2008rigid}. Also, in the context of network localization for mobile sensors, rigidity is a basic condition for successful localization \cite{aspnes2006theory}.   Furthermore, as  will be shown in later analysis, connectivity maintenance and collision avoidance can also be achieved as direct consequences of the rigidity preservation for mobile robot networks.

   There exist several  fundamental theorems relevant for  rigidity testing, and in  this paper we focus on those using linear algebra  for a quantitative rigidity analysis provided by the rigidity matrix. Inspired by the concept of the Fiedler (algebraic connectivity) eigenvalue of the Laplacian  matrix, we also adopt a particular eigenvalue characterization of the rigidity level for a given  framework.   To this end, an extended version of the Laplacian matrix which we term the {\it{rigidity Laplacian matrix}} is constructed  by incorporating  robots' positions and the network topology of the formation.

The idea of using the eigenvalue information to describe the rigidity property  has something in common with the  quantitative measure  using stiffness matrix \cite{kim2010optimizing}, \cite{zhu2011link} and an anchor selection metric  from  a reduced-order rigidity matrix in the localization optimization problem \cite{shames2009minimization}. A more recent effort towards this direction which is closer to the idea of this paper can be found in \cite{zelazo2013decentralized}. However, the analysis and results in this paper are distinguished from these previous results in several ways. First, we derive a simple and general form of the rigidity matrix which involves both positional configuration and network topology, and thus the analysis on the eigenstructure of the  rigidity Laplacian matrix can reveal much interesting information for the global formation. Furthermore,  some other useful properties of the rigidity index are also provided, which have not been covered elsewhere.
Another contribution of this paper  is the development of  an entirely distributed estimation and control solution for the rigidity preservation. This can be seen as  parallel to work  dealing with distributed connectivity  control \cite{yang2010decentralized}, \cite{sabattini2012decentralized}, \cite{sabattin2013}. As an extra novelty,  this paper also highlights the design of a distributed estimation algorithm where the convergence rate is controllable and is distinctly faster than the distributed power iteration scheme in other papers \cite{yang2010decentralized}, \cite{sabattini2012decentralized}, \cite{sabattin2013}.

This paper is organized as follows. Section II reviews some graph theoretic preliminaries as well as basic concepts on formation rigidity. In Section III, by deriving a simple form for the rigidity matrix, we discuss some interesting properties of the  rigidity Laplacian matrix. The problem of rigidity preservation is formulated as to guarantee the positivity of the desired eigenvalue. Section IV discusses a potential function-based framework for this rigidity preservation problem. A distributed estimation solution via an inverse power iteration method is discussed in detail in Section V. The extension to the 3-D case is briefly examined in Section VI. In Section VII simulation results are provided, and concluding remarks are provided in Section VIII.

\subsection{Notations}
The notations used in this paper are fairly
standard. $\mathbb{R}^n$ denotes the $n$-dimensional Euclidean space. $\mathbb{R}^{m\times n}$
denotes the set of $m\times n$ real matrices. If $M$ is
a vector or matrix, its transpose is denoted by $M^T$. The rank, image and null space of matrix $M$ are denoted by  $rank(M)$, $Im(M)$ and $null(M)$, respectively. For a symmetric matrix $M$, its $i$-th smallest eigenvalue is denoted by $\lambda_i(M)$. The notation $\text{diag}\{x\}$ denotes a (block) diagonal matrix with the (block) vector $x$ on its diagonal. $span\{v_1, v_2, \cdots, v_k\}$ represents the subspace spanned by a set of vectors $v_1, v_2, \cdots, v_k$. $I_n$ is the $n \times n$ identity matrix, and  $\mathbf{1}_n$ denotes a $n$-tuple column vector of all ones.   The symbol $\otimes$ denotes the Kronecker product.

\section{Preliminaries}
In this section we introduce some basic notations and concepts on graph theory and rigidity theory. Further details can be found in \cite{mesbahi2010graph} and \cite{anderson2008rigid}.

\subsection{Graph theory}
We assume that the mobile robots are modeled by kinematic points. Consider an undirected graph with $m$ edges and $n$ vertices, denoted by $\mathcal{G} =( \mathcal{V, \mathcal{E}})$  with vertex set $\mathcal{V} = \{1,2,\cdots, n\}$ and edge set $\mathcal{E} \subset \mathcal{V} \times \mathcal{V}$.  The vertex set represents the robots (and we may use the word \emph{agent} interchangeably in the context) and the edge set represents the communication links between different robots. The matrix relating the nodes to the edges is called the incidence matrix $H = \{h_{ij}\} \in \mathbb{R}^{m \times n}$, whose entries are defined as (with arbitrary edge orientations)
     \begin{equation}
     h_{ij} =  \left\{
       \begin{array}{cc}
       1,  &\text{ the } i\text{-th edge sinks at node }j  \\ \nonumber
       -1,  &\text{ the } i\text{-th edge leaves at node }j  \\ \nonumber
       0,  & \text{otherwise}  \\
       \end{array}
      \right.
      \end{equation}

In the above definition we use the convention that each row of $H$ represents one existing edge which links two vertices.

The adjacency matrix $A(\mathcal{G})$ is a symmetric $n\times n$ matrix encoding  the vertex adjacency relationships, with entries $A_{ij} = 1$ if $\{i,j\} \in \mathcal{E}$, and $A_{ij} = 0$ otherwise.   Another important matrix representation of a graph $\mathcal{G}$ is the Laplacian matrix $L(\mathcal{G})$, which is defined as $L(\mathcal{G}) = H^TH = \text{diag}\{A\mathbf{1}\}-A$.

Some properties of the graph Laplacian matrix $L(\mathcal{G})$ are summarized in the lemma below \cite{mesbahi2010graph}.
\begin{lemma}
 Given an undirected graph $\mathcal{G}$:
\begin{itemize}
\item $L(\mathcal{G})$  is orientation-independent.
\item $L(\mathcal{G})$ is symmetric and positive semidefinite.
\item If $\mathcal{G}$ is connected, then $L(\mathcal{G})$ has one and only one zero eigenvalue, with $null(L(\mathcal{G})) = span\{\mathbf{1}\}$.
\item  $ x^TLx = \sum_{\{i,j\} \in \mathcal{E}} A_{ij} (x_i-x_j)^2$ where $x$ is a column vector.
\end{itemize}
\end{lemma}

\subsection{Formation graph and infinitesimal rigidity}
We embed the graph $\mathcal{G}$ into  2-D space (in this paper we mostly focus on the analysis in 2-D space, however the extension to 3-D space is straightforward). Let $p_i = [p_{ix},\, p_{iy}]^T \in \mathbb{R}^2, i \in \{1,2,\cdots, n\}$ denote the position of node $i$. The stacked  vector $p=[p_1^T, \,
p_2^T, \cdots, \, p_n^T]^T$ represents the position configuration for all the  $n$ nodes. By introducing the matrix $\bar H: = H\otimes I_2 \in \mathbb{R}^{2m\times2n}$, one can construct the edge space as an image of $\bar H$ from the position vector $p$:
\begin{equation}
z = \bar H p  \label{z_equation}
\end{equation}
with $z_i = [z_{i,x}, z_{i,y}]^T \in \mathbb{R}^2$ being the relative position vector for the vertex pair defined by the $i$-th edge. In the following, two notations, $z_k$ and $z_{k_{ij}}$ will be used interchangeably to denote the $k$-th edge which links agent $i$ and agent $j$.

The rigidity function $r_{\mathcal{G}}(p): \mathbb{R}^{2n} \rightarrow \mathbb{R}^m$ associated with the framework $(\mathcal{G}, p)$ is defined as:
\begin{equation}
r_{\mathcal{G}}(p) = \frac{1}{2} \left[\cdots, \|p_i-p_j\|^2,  \cdots \right]^T
\end{equation}
where the norm is the standard Euclidean norm, and the $k$-th component in $r_{\mathcal{G}}(p)$, $\|p_i-p_j\|^2$, corresponds to the square length of edge $z_k$. The framework $(\mathcal{G}, p)$ is said to be rigid, if there exists an open neighbourhood $\mathcal{U}$ of $p$ such that, if $q \in \mathcal{U}$ and $r_{\mathcal{G}}(p) = r_{\mathcal{G}}(q)$, then $(\mathcal{G}, p)$ is congruent to $(\mathcal{G}, q)$.

Another useful tool to characterize the rigidity property of a framework is the rigidity matrix, which is defined as
\begin{equation}
R(p) = \frac{\partial r_{\mathcal{G}}(p) } {\partial p}  \label{R_matrix}
\end{equation}

 The framework $(\mathcal{G}, p)$ is said to be infinitesimally rigid if the rank of the rigidity matrix $R$ equals $2n-3$. Also, if $(\mathcal{G}, p)$ is  infinitesimally rigid, so is $(\mathcal{G}, p')$ for a generic (open and dense) set of $p'$. Generally speaking, infinitesimal rigidity implies rigidity, but the converse is not true. In the rest of this paper, we will use the rank condition of the rigidity matrix to determine whether a formation is rigid.

\subsection{Communication model with limited sensing range}
Suppose that agents $i$ and $j$ are able to interact with each other  if their distance is within a communication radius $\kappa$. The communication topology can then be modeled by an undirected \emph{dynamic} graph $\mathcal{G} = \{\mathcal{V}, \mathcal{E}\}$, where $\mathcal{E} \subseteq \mathcal{V} \times \mathcal{V}$ denotes the set of communication links:
\begin{equation}
\mathcal{E} = \mathcal{E}(\mathcal{G}) = \{\{i,j\}| \|p_i - p_j\| \leq \kappa,  i, j \in \mathcal{V}, i \neq j\}
\end{equation}

Thus, the neighbors of the $i$-th agent are given by
\begin{equation}
\mathcal{N}_i = \{j \in \mathcal{V} | \{i,j\} \in \mathcal{E}\}
\end{equation}

It is also desirable to define weighted edges for the graph. The weight should be a function of the distance
between agent pairs. Some choices and discussion for different weight functions can be found in \cite{kim2006maximizing}. We consider the following edge weight function
     \begin{equation}
     w_{k_{ij}} =  \left\{
       \begin{array}{cc}
       e^{-\|z_{k_{ij}}\|^2/(2\sigma^2)},  &\,\,\text{if } \|z_{k_{ij}}\| \leq  \kappa  \\
       0,  & \,\,\, \text{otherwise}
       \end{array}
      \right.
      \end{equation}
      The weight will decrease when the inter-agent distance gets larger. One can choose the scalar parameter $\sigma$ to satisfy a threshold condition $e^{-\kappa^2/(2\sigma^2)} = \sigma' $, with $\sigma'$ being a small predefined threshold.

By defining an $m\times m$  diagonal matrix $W$ whose diagonal entries are the weights for each edge,  one can construct the weighted Laplacian matrix $L(\mathcal{G}_{w}) = H^TWH$. All the properties stated in Lemma 1 also apply to the weighted Laplacian matrix $L(\mathcal{G}_{w})$ \cite{mesbahi2010graph}.

\section{Rigidity Laplacian matrix and rigidity preservation index}

\subsection{Rigidity Laplacian matrix}
Firstly we would like to derive a simple  expression for the rigidity matrix which involves both the network topology and position configuration. Recall  \eqref{z_equation}, which shows that the edge space lies in the image of $\bar H$. The rigidity function is a map from the node positions to the squared edge lengths. Thus we can redefine the rigidity function, $g_{\mathcal{G}}(z): Im(\bar H) \rightarrow \mathbb{R}^m$ as $g_{\mathcal{G}}(z) = \frac{1}{2}\left[\|z_1\|^2, \|z_2\|^2, \|z_3\|^2, \cdots, \|z_m\|^2 \right]^T$. From \eqref{z_equation} and \eqref{R_matrix}, one can obtain the following simple form for the rigidity matrix
\begin{eqnarray} \nonumber
R(p) = \frac{\partial r_{\mathcal{G}}(p)}{\partial p} &=& \frac{\partial g_{\mathcal{G}}(z)}{\partial z} \frac{\partial z} {\partial p}  \\ \nonumber
&=& \left(\begin{array}{ccc}
z_{1}^T & \cdots & 0  \\
\vdots & \ddots & \vdots  \\
0 & \cdots & z_{m}^T   \\
\end{array}\right)
\bar H  \\
&=& Z^T \bar H
\end{eqnarray}
where $Z$ is  a block diagonal matrix $Z = \textrm{diag} \{z_1,\,z_2,\cdots,\,z_m\}$.

The set of all infinitesimal displacements caused by the  rigid body motions forms a subspace of dimension three, which also serves as the null space of the rigidity matrix $R$. In fact, a set of linearly independent null vectors of  $R$ can be calculated directly as
\begin{eqnarray}
v_1 &=& \mathbf{1}\otimes [1, 0]^T = [1, 0, 1,0, \cdots, 1,0]^T \\  \label{null_1}
v_2 &=& \mathbf{1}\otimes [0, 1]^T = [0, 1,0, 1, \cdots, 0,1]^T  \\  \label{null_2}
v_3 &=& [p_{1y}, -p_{1x}, p_{2y}, -p_{2x}, \cdots, p_{ny}, -p_{nx}]^T  \label{null_3}
\end{eqnarray}

Following the definition of the weighted Laplacian matrix of $L(\mathcal{G}) = H^TWH$, we construct a new matrix $E$ in a similar way: $E(\mathcal{G}, p) = R^TWR = \bar H^T ZWZ^T\bar H$. Since the matrix $E$ shares several similar properties with  $L$, we  term it  the \emph{rigidity Laplacian matrix}. In fact,  the matrix $E$ can be regarded as a position-weighted Laplacian matrix for the  framework $(\mathcal{G}, p)$, where the weights are described by a diagonal block matrix $ZWZ^T$ involving the position information. The block diagonal matrix $ZWZ^T$ is expressed by
\begin{equation}
ZW Z^T = \text{diag}\{w_1z_1z_1^T,w_2z_2z_2^T,\cdots, w_mz_mz_m^T\}
\end{equation}
where $w_i$ is a scalar and $z_iz_i^T$ is a $2 \times 2$ block:　
\begin{equation} z_iz_i^T =
\left(\begin{array}{cc}
z_{i,x}^2 & z_{i,x}z_{i,y} \\
z_{i,x}z_{i,y} & z_{i,y}^2  \\
\end{array}\right)
\end{equation}

The structure of $E$ also resembles that of the Laplacian matrix $L$, while
$E$ has twice the dimension compared to $L$. In fact, it is more convenient to consider the $2\times 2$ block entries of  matrix $E = [E_{ij}]_{1\leq i,j \leq n} \in \mathbb{R}^{2n \times 2n}$:
     \begin{equation}
     E_{ij} =  \left\{
       \begin{array}{cc}
       \sum_{l \in \mathcal{N}_i} w_{il}(z_{k_{il}}z_{k_{il}}^T), & \,\,\, \text{if} \,\, i=j  \\
        - w_{ij}(z_{k_{ij}}z_{k_{ij}}^T), & \,\,\, \text{if}\,\, i\neq j \,\, \text{and} \,\,\{i,j\} \in \mathcal {E}  \\
        \mathbf{0}_{2 \times 2}, & \,\,\, \text{if}\,\,  i\neq j  \,\, \text{and} \,\,\{i,j\} \notin \mathcal {E} \\
       \end{array}  \label{struc_E}
      \right.
      \end{equation}

\parskip= 8pt
The following result shows that the rigidity Laplacian matrix $E$ shares the same null space with rigidity matrix $R$.
\begin{lemma}
$null(E(\mathcal{G}, p)) = null(R^TR) = null(R)$.
\end{lemma}

\noindent \begin{proof}
It is obvious that $null(R^TR) = null(R)$. The weight matrix $W$ is invertible as it is a diagonal matrix with positive diagonal entries. Thus $rank(E(\mathcal{G}, p)) = rank(R^TR)$ and the eigenspace corresponding to the zero eigenvalues of $E$ is the same as that of $R^TR$.
\end{proof}
The above result will be used in Section V to construct a modified matrix based on the null vectors of $R$ for the distributed eigenvector estimation of $E$.

Similarly to Lemma 1 on properties of Laplacian matrix $L$,  some properties of $E(\mathcal{G})$ are listed in the following lemma.
\begin{lemma}
Given a framework $(\mathcal{G}, p)$:
\begin{itemize}
\item $E(\mathcal{G}, p)$ is orientation-independent.
\item $E(\mathcal{G}, p)$ is symmetric and positive semidefinite.
\item If the framework $(\mathcal{G}, p)$ is infinitesimally rigid, then $E(\mathcal{G}, p)$ has three and only three zero eigenvalues, with $null (E) = span\{v_1,  v_2, v_3\}$.
\end{itemize}
\end{lemma}

\subsection{Rigidity preservation index}

The second smallest eigenvalue of the Laplacian matrix, $\lambda_2(L)$, also called the algebraic connectivity value or the Fiedler eigenvalue, plays an important role for network analysis. As a straightforward extension, we choose the critical eigenvalue $\lambda_4(E)$ as a quantitative
index of the rigidity level for the formation $(\mathcal{G}, p)$. This idea is similar to \cite{kim2010optimizing}, \cite{zhu2011link}, where a worst rigidity index  for  the  stiffness matrix was defined. In a recent work  \cite{zelazo2013decentralized} this eigenvalue  was also used for measuring the rigidity property of the embedded framework.

The control problem of preserving network connectivity, which is to guarantee $rank(L) = n-1$ or $\lambda_2(L)>0$, has been extensively studied in the literature \cite{zavlanos2011graph}. The rigidity preservation problem can be formulated in a similar way: designing control schemes to ensure that $rank(E) = rank(R) = 2n-3$, or equivalently, to guarantee that $\lambda_4(E)>0$ for a group of mobile robots moving in an uncertain environment.

As a consequence of preserving the rigidity for mobile robot teams, some other nice formation properties can also be achieved, which are summarized in the following propositions.

\parskip= 8pt
\begin{proposition}
$\lambda_4(E)>0$  implies that $\lambda_2(L)>0$. That is,  the maintenance of rigidity implies the maintenance of graph connectivity.
\end{proposition}

\noindent \begin{proof}
Since $\lambda_4(E)>0$, then $rank(E) = 2n-3$.  From the definition $E=\bar{H}^TZWZ^T\bar{H}$, one has $rank(\bar{H}) \geq rank(E) = 2n-3$. It follows that $rank(H)\geq n-1.5$. Since the maximum rank of $H$ is $n-1$, this immediately implies that $rank(H) = n-1$.  Therefore the graph is connected.
\end{proof}

\parskip= 8pt
\begin{proposition}
$\lambda_4(E)>0$ also implies collision avoidance between each pair of neighbor agents.
\end{proposition}

\noindent \begin{proof}
Suppose neighboring agent  $i$ and agent $j$ collide, then one has $z_{k_{ij}} = 0$. This will introduce a zero block in the matrix $ZWZ^T$ and thus dim$(null(E))>3$, which violates the condition of $\lambda_4(E)>0$. Thus, the constraint of $\lambda_4(E)>0$ implies that no $z_{k_{ij}}$ will be zero, i.e. no agent pairs will collide.
\end{proof}

\section{Energy function based control approach}
Consider  a group of $n$ robots whose dynamics are described by the single-integrator model
\begin{equation}
\dot p_i = u_i^r
\end{equation}
where $u_i^r$ is the designed control input for the $i$-th robot. As stated in Section III, the rigidity preservation is equivalent to guaranteeing that $\lambda_4(E)$ is strictly greater than zero.   Let $\epsilon$ be a desired lower bound for $\lambda_4(E)$. The control objective is then to ensure that $\lambda_4$ never goes below this lower bound. Inspired by \cite{zavlanos2007potential}, we also employ the energy function method to construct potential fields for generating decentralized control strategies. Denote $V(\lambda_4): \mathbb{R}^+ \rightarrow \mathbb{R}^+ $ as a positive definite  energy  function of $\lambda_4$ over the interval $(\epsilon, +\infty)$. The energy function $V$ is defined such that the following properties hold:
\begin{itemize}
\item It is continuously differentiable.
\item It is non-negative.
\item It is non-increasing with respect to $\lambda_4$.
\item When $\lambda_4 \rightarrow \epsilon$, $V \rightarrow + \infty$ and $\|\frac{\partial V}{\partial \lambda_4}\| \rightarrow + \infty$.
\item When $\lambda_4 \rightarrow + \infty$, $V$ approaches a constant with a vanishing slope.
\end{itemize}
A good choice for the energy function can be
\begin{equation}
V(\lambda_4) = coth(\lambda_4 - \epsilon)
\end{equation}
where $coth$ is the Hyperbolic Cotangent function. This energy function has also been  used for the network algebraic connectivity control \cite{sabattini2012decentralized}, \cite{sabattin2013}, \cite{ji2007distributed}, \cite{sabattini2011distributed}.

The control design essentially drives the robot teams to perform a gradient descent of $V(\cdot)$ to maximize the value of $\lambda_4$ such that the rigidity preservation can be guaranteed:
\begin{equation}
u_i^r =  -\frac{\partial V} {\partial p_i} = -\frac{\partial V} {\partial \lambda_4} \frac{\partial \lambda_4} {\partial p_i}
\end{equation}

 Denote  $v_4 = [v_{4,1}^T, v_{4,2}^T, \cdots, v_{4,n}^T] \in \mathbb{R}^{2n}$ with $v_{4,i} = [v_{4,ix},v_{4,iy}]^T$ as a normalized eigenvector corresponding to the eigenvalue $\lambda_4$ of the matrix $E$. One has $\lambda_4(E) = v_4^TEv_4$.
In the following, we would like to obtain the closed form of the term $\partial \lambda_4 / \partial p_i$.  According to the structure of the matrix $E$ in \eqref{struc_E}, one has
\begin{eqnarray}
v_4^TEv_4 &=& \sum_{i,j} v_{4,i}^T E_{ij} v_{4,j} \\ \nonumber
&=&\sum_{\{i,j\} \in \mathcal{E} } w_{ij} z_{ij}^T (v_{4,i} - v_{4,j}) (v_{4,i} - v_{4,j})^T z_{ij}\\ \nonumber
\end{eqnarray}

Hence, the explicit expression should be
\begin{small}
 \begin{eqnarray}
\frac{\partial (v_4^TEv_4)} {\partial p_i} &=& v_4^T\frac{\partial E}{\partial p_i}v_4 \\ \nonumber \label{control}
  &=&  2\sum_{j \in \mathcal{N}_i} w_{ij} \frac{\partial z_{ij}^T}{\partial p_i} (v_{4,i} - v_{4,j}) (v_{4,i} - v_{4,j})^T z_{k_{ij}} \\ \nonumber
  && + \sum_{j \in \mathcal{N}_i} \frac{\partial w_{ij}}{\partial p_i} z_{k_{ij}}^T (v_{4,i} - v_{4,j}) (v_{4,i} - v_{4,j})^T z_{k_{ij}} \\ \nonumber
    &=&  2\sum_{j \in \mathcal{N}_i} w_{ij} (v_{4,i} - v_{4,j}) (v_{4,i} - v_{4,j})^T (p_i-p_j) \\ \nonumber
    &&+  \sum_{j \in \mathcal{N}_i} \frac{\partial w_{ij}}{\partial p_i} z_{k_{ij}}^T (v_{4,i} - v_{4,j}) (v_{4,i} - v_{4,j})^T z_{k_{ij}} \\ \nonumber
\end{eqnarray}
\end{small}
In the first line we have used the fact that $\frac{\partial (v_4^TEv_4)} {\partial p_i} = v_4^T\frac{\partial E}{\partial p_i}v_4$ due to $\|v_4\|=1$ and $E = E^T$ \cite{parlett1980symmetric}.

The final control for each agent should also contain some additional terms for achieving extra tasks such as  flocking, moving or  searching  algorithms, etc. Though the above gradient control is designed for the aim of distributed implementation, one can find that the control input  still contains some global information such as  $v_4$ and $\lambda_4$.  This may require a centralized controller to collect the information from all the agents, do the calculation, and then broadcast it to all the agents. However, for  large-scale robot networks, it is not desirable to implement such a centralized controller. This motivates us to find some distributed estimation algorithm and to design an entirely distributed control scheme.

 \section{Distributed estimation using inverse power iteration}
 A key problem in designing  a distributed  control is to estimate the global information $v_4$ and $\lambda_4$ for each robot in a local way. The power iteration method is an established method for estimating the dominant eigenvector for a specific matrix, assuming there is a single eigenvalue of
maximum modulus \cite{parlett1980symmetric}, \cite{yang2010decentralized}. The shifted inverse iteration method, which is a variation  of the power iteration method,  can be used to estimate any eigenvalue (instead of the dominant eigenvalue), provided that a suitably accurate initial estimate of the desired eigenvalue is given. If the estimate is very close to the desired eigenvalue, this inverse power iteration method is generally much faster than the standard power iteration method.

The reason for choosing the inverse power iteration method is to improve the convergence speed of the estimation process. Since all the robots are in a dynamic environment and their positions and communication
links are under change from time to time, it is desirable to devise a fast iteration estimation scheme  to satisfy the control requirements. Also for this reason, the inverse power iteration has been discussed  in a recent work
\cite{williams2013locally} for the constrained connectivity control.

Before presenting the estimation procedure via the inverse power iteration method, we firstly introduce two closely-related distributed algorithms  which will be incorporated in the estimation procedure of the desired eigenvector.

\textbf{\emph{Distributed average consensus estimator}}\\
The following dynamic proportional-integral consensus estimator will be used frequently for average value calculation \cite{freeman2006stability}, \cite{yang2010decentralized}:
\begin{eqnarray*}
\dot z^i &=& \rho (\alpha^i - z^i) - K_P \sum_{j \in \mathcal{N}_i}(z^i - z^j) + K_I \sum_{j \in \mathcal{N}_i}(w^i - w^j) \\ \nonumber
\dot w^i &=& -K_I \sum_{j \in \mathcal{N}_i}(z^i - z^j)
\end{eqnarray*}
where $\alpha^i$ is some time-varying measurement, $z_i$ is  an estimate at node $i$ of the average value of the $\alpha^j$ over all nodes, $\rho > 0$ is the rate new information replaces old information and $K_P$, $K_I$ are estimator gains. As shown in \cite{yang2010decentralized}, the gains should be chosen large enough such that the time constant for consensus estimation is much less than the time constants for eigenvector estimation and motion controllers. This consensus estimator allows $n$ agents to compute an average approximation $\bar z = \frac{1}{n}\sum_{i} \alpha^i$ by using only local interaction even if $\alpha^i$ are varying. In the following, we use the function symbol $\text{Ave}(\cdot)$ to denote the average consensus operation using this estimator.

\textbf{\emph{Jacobi overrelaxation method}}\\
The inverse power iteration method involves the inverse calculation of a specific matrix, which can be transformed to the problem of  solving a linear equation.  Consider  the following linear equation:
\begin{equation}
Ax = b  \label{linea}
\end{equation}

There are several numerical algorithms available for solving \eqref{linea} in a parallel and distributed way \cite{bertsekas1997parallel}. One of the powerful algorithms is called the Jacobi overrelaxation method, which involves the following  iterative steps:
\begin{equation}
x_i(t+1) = (1-\gamma)x_i(t)-\frac{\gamma}{a_{ii}}\left( \sum_{j\neq i}^{n}a_{ij}x_{j}(t)-b_i \right)
\end{equation}
If A is symmetric and positive definite and $\gamma>0$ is sufficiently small, then the sequence generated by the above algorithm converges to a solution of  $Ax = b$  \cite{bertsekas1997parallel}.

\subsection{Distributed eigenvector estimation using shifted inverse power iteration}

Let us consider the matrix $\bar{E}: = E + Q$ where $Q = \vartheta v_1 v_1^T + \vartheta v_2v_2^T + \vartheta v_3v_3^T$ and $\vartheta$ is some sufficiently large positive constant. Note that $Q$ is a rank-three symmetric matrix. It has three positive eigenvalues and the rest are all zero. Denote the positive eigenvalues of $Q$ as $\eta_1, \eta_2, \eta_3$. The value of $\vartheta$ is chosen such that $min\{\eta_1, \eta_2, \eta_3\} > \lambda_{4}$. This can be determined  by doing some prior calculation in advance using the value $\lambda_4^{max}$.  The eigenstructure of the matrix $\bar E$ is as follows (note that the list is not necessarily in an ascending order):
\begin{equation}
\text{Eigenvalues of } \bar{E}: \lambda_{4},\,\,  \lambda_{5},\cdots, \lambda_{2n},\,\,  \eta_1, \,\, \eta_2,\,\,  \eta_3
\end{equation}
with associated eigenvectors (or eigenspace)
\begin{equation}
 v_{4}, \,\,  v_{5}, \cdots, v_{2n}, \,\, span\{ v_{1}, \,\, v_{2},\,\,  v_3\}
\end{equation}

Define a new matrix $(\bar{E} - \mu I_{2n})^{-1}$ where $\mu$ is a positive number. It follows that the set of  eigenvectors of $(\bar{E} - \mu I_{2n})^{-1}$ are the same as those of $\bar E$, with the
eigenvalues listed below:
\begin{equation}
(\lambda_{4}-\mu)^{-1}, \cdots, (\lambda_{2n} -\mu)^{-1}, (\eta_1 - \mu)^{-1}, (\eta_2 - \mu)^{-1}, (\eta_3 - \mu)^{-1}
\end{equation}

By choosing $\mu$ close to $\lambda_{4}$, the dominant eigenvalue of $(\bar{E} - \mu I_{2n})^{-1}$ will be much larger in magnitude than other eigenvalues listed above. Thus, by doing the
power iteration method on the  matrix $(\bar{E} - \mu I_{2n})^{-1}$, the convergence rate can be greatly improved.

Another issue is to obtain in a distributed way the inverse of the matrix $(\bar{E} - \mu I_{2n})$. Instead of doing the matrix inversion operation, we would like to solve the following linear equation:
\begin{equation}
(\bar{E} - \mu I_{2n}) r = {\widetilde {v}_4}^{(k)}  \label{linear_eq}
\end{equation}
where ${\widetilde {v}_4}^{(k)}$ is the estimate of $v_4$ at the $k$-th step.
 The Jacobi overrelaxation iteration method is employed
to solve it in a distributed way:
\begin{small}
\begin{eqnarray}  \label{iteration}
r_i^{(p+1)} &=& (1-\gamma)r_i^{(p)} \\ \nonumber
 &&- \frac{\gamma}{\bar{E}_{ii}-\mu} \left(- {\widetilde {v}}_{4,i}^{(k)}-\sum_{j \in \mathcal{N}_i} E_{ij} r_{j}^{(p)} + \sum_{j\neq i}Q_{ij} r_j^{(p)} \right)
\end{eqnarray}
\end{small}
where $Q_{ij}$ is the $ij$-th block of the matrix $Q$:
\begin{equation} Q_{ij} = \vartheta
\left(\begin{array}{cc}
p_{iy}p_{jy}+1  & -p_{iy}p_{jx} \\
-p_{ix}p_{jy} & p_{ix}p_{jx}+1  \\
\end{array}\right)
\end{equation}
The iteration step in \eqref{iteration} involves mostly local communication except that the calculation of the last term $\sum_{j\neq i}Q_{ij} r_j^{(p)}$ requires global information. Note that $\sum_{j\neq i}Q_{ij} r_j^{(p)} = n \text{Ave}(\sum_{j} Q_{ij} r_j^{(p)}) - Q_{ii}r_i^p $, and $\text{Ave}(\sum_{j} Q_{ij} r_j^{(p)})$ can be computed by using the local average estimator
\begin{small}
\begin{equation}   \vartheta
\left(\begin{array}{cc}
p_{iy} \text{Ave}(p_{iy}r_{i,x}^{(p)})+\text{Ave}(r_{i,x}^{(p)})   -p_{iy}\text{Ave}(p_{ix}r_{i,y}^{(p)}) \\  \label{ave_Q}
-p_{ix}\text{Ave}(p_{iy}r_{i,x}^{(p)}) + p_{ix}\text{Ave}(p_{ix} r_{i,y}^{(p)})+\text{Ave}(r_{i,y}^{(p)})  \\
\end{array}\right)
\end{equation}
\end{small}
Suppose after $\bar p$ steps the iteration of the solution to \eqref{linear_eq} converges.
The followed normalization step can also be implemented in a distributed way by using the average consensus estimator:
\begin{equation}
{\widetilde {v}}_{4,i}^{(k+1)} = \frac{r_i^{(\bar p)}}{\sqrt {{n\text{Ave}((r_{i,x}^{(\bar p)})^2}+(r_{i,y}^{(\bar p)})^2) }}
\end{equation}
Further suppose that after $\bar k$ steps, the convergence of the eigenvector estimation is achieved. Then the eigenvalue can be estimated by using the Rayleigh quotient
\begin{equation}
{\widetilde {\lambda}}_4 = ({\widetilde {v}_4}^{(\bar k)})^T E ({\widetilde {v}_4}^{(\bar k)})
\end{equation}

However, the actual value ${\widetilde {\lambda}}_4$ cannot be computed by each agent as the information of the estimated normalized eigenvector ${\widetilde {v}}^{(\bar k)}$ cannot be accessed by all the agents. Nevertheless, the local estimation of the eigenvalue can still be computed by using again the average consensus procedure. The initial input for the average consensus is
\begin{equation}
z_i(0) =  ({\widetilde {v}}_{4,i}^{(\bar k)})^T \sum_{j \in \mathcal{N}_i} E_{ij}({\widetilde {v}}_{4,i}^{(\bar k)} - {\widetilde {v}}_{4,j}^{(\bar k)})
\end{equation}
Thus robot $i$ can calculate its local estimation of the eigenvalue by
\begin{equation}
{\widetilde {\lambda}}_{4,i} = n \text{Ave}(z_i)
\end{equation}
Then the control input in (18) can be modified by replacing $\lambda_4$ and $v_4$ with their estimates.

\subsection{Conditions and Convergence rate of the estimation}
Two important aspects of the estimation procedure should be emphasized.
When the Jacobi overrelaxation  iteration method is used to solve the linear equation related to the matrix $(\bar{E} - \mu I_{2n})$, the condition for  convergence is that $\gamma>0$ is sufficiently small, and the matrix $(\bar{E} - \mu I_{2n})$ should be symmetric and positive definite (Page 154 of \cite{bertsekas1997parallel}). The parameter $\gamma$ can be adjusted in the implementation process. In order to ensure the positive definiteness of $(\bar{E} - \mu I_{2n})$, we can choose $\vartheta$ and $\mu$ so that $min\{\eta_1, \eta_2, \eta_3\} > \lambda_{4}$ and $0<\mu<\lambda_4$. The latter condition is to ensure that $(\lambda_4 - \mu)^{-1}$ is the dominant eigenvalue of the matrix $(\bar{E} - \mu I_{2n})^{-1}$. In practise the true value for $\lambda_4$ is unknown, but a lower bound for $\lambda_4$ is known to be $\varepsilon$. Thus, a conservative range for the shift $\mu$ can always be chosen as $0<\mu<\varepsilon$.

The convergence rate of the estimation of the eigenvector is controlled by the ratio $|\frac{\lambda_4 - \mu}{\lambda_5 - \mu}|$. Thus, by choosing  $\mu$  closer to $\lambda_4$, the convergence rate will be faster.
In the implementation, the  initial value $\mu$ can be chosen by combining  other methods. For example, by employing a few steps of power iteration, some reasonable initial guess of $\lambda_4$ can be obtained.

\section{Extension to 3-D case}

The problem description and notations in the 3-D case are similar to those in above sections. The matrix $E$ is modified as $E = R^TR = {\bar{H}}^T ZW Z^T{\bar{H}}$, where $\bar{H} = H
\otimes I_{3}$.  The block diagonal matrix $Z WZ^T$ is expressed as
\begin{equation}
Z Z^T = \text{diag}\{w_1 z_1z_1^T, w_2 z_2z_2^T,\cdots, w_m z_mz_m^T\}
\end{equation}
where $z_iz_i^T$ is a $3 \times 3$ block:　
\begin{equation} z_iz_i^T =
\left(\begin{array}{ccc}
z_{i,x}^2 & z_{i,x}z_{i,y} & z_{i,x}z_{i,z} \\
z_{i,x}z_{i,y} & z_{i,y}^2 & z_{i,y}z_{i,z} \\
z_{i,x}z_{i,z} & z_{i,y}z_{i,z} & z_{i,z}^2\\
\end{array}\right)
\end{equation}

The null vectors for the matrix $R$ in the 3-D space case are listed as below:
\begin{eqnarray}
v_1 = \mathbf{1}\otimes [1, 0, 0]^T = [1, 0,0, 1,0,0, \cdots, 1,0,0]^T   \\
v_2 = \mathbf{1}\otimes [0, 1, 0]^T = [0,1, 0,0, 1,0, \cdots, 0,1,0]^T \\
v_3 = \mathbf{1}\otimes [0, 0, 1]^T = [0,0,1, 0, 0, 1, \cdots,0,0, 1]^T \\
v_4 = [p_{1y}, -p_{1x}, 0, p_{2y}, -p_{2x}, 0, \cdots, p_{ny}, -p_{nx}, 0]^T \\
v_5 = [p_{1z},  0, -p_{1x}, p_{2z}, 0, -p_{2x},  \cdots, p_{nz}, 0, -p_{nx}]^T \\
v_6 = [0, p_{1z}, -p_{1y}, 0, p_{2z}, -p_{2y},  \cdots, 0, p_{nz}, -p_{ny}]^T
\end{eqnarray}

Note that if the rank of the rigidity matrix satisfies $rank(R) = 3n-6$, then the framework imbedded in the 3-D space is infinitesimally rigid. Thus, the critical eigenvalue for the 3-D case is $\lambda_7$. All the analysis  above  can be applied to the 3-D case, with only slight modifications required.

\section{Simulation results}
\subsection{Comparisons between distributed power iteration method and  inverse power iteration method}
In this section we give some simulation results via Matlab/Simulink. Suppose we have a system of $n=5$ robots operating in a bounded workspace in the plane. The communication radius is $\kappa = 10$, the threshold is set as $\sigma' =0.01$ and the relaxation parameter is chosen $\gamma = 0.25$. Without loss of generality, we simulate using randomly-generated  positions for the robots in a $10\times 10$ square area, which generates a matrix $E$ with the spectrum $[0, 0, 0, 4.46, 8.49, 11.18, 18.56, 22.85,  32.74, 38.93]$.
By using the distributed inverse power iteration method, the desired eigenvector has been estimated by each agent (simulation results not shown here).  We consider the power iteration method discussed in \cite{yang2010decentralized}, \cite{sabattin2013}, \cite{zelazo2013decentralized} and modify a discrete-time version for the comparison. To compare the convergence speed, all the assumptions and initial conditions of the estimates are the same. The results are shown in Fig.1. It is obvious that the distributed inverse power iteration method proposed in this paper displays superior performance over the distributed power iteration method proposed in \cite{yang2010decentralized}, \cite{sabattin2013}, \cite{zelazo2013decentralized}.

Also from Fig.1 one can observe that, when the initial guess $\mu$ is chosen closer to the true value of $\lambda_4$, the convergence will be much faster, which can be achieved by only a few iteration steps. This property is quite favorable for distributed large-scale robot network control.

\begin{figure}
  \centering
  \includegraphics[width=85mm]{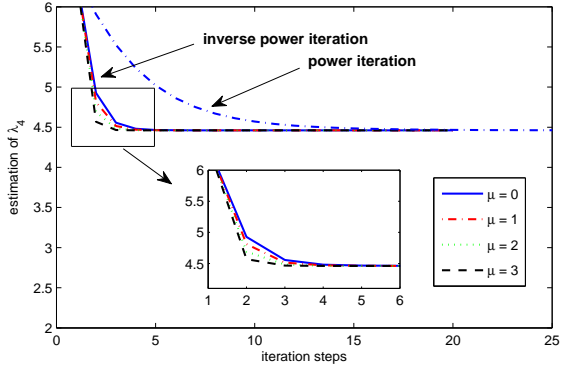}
  \caption{Convergence speed of the eigenvalue estimation between power iteration and inverse power iteration method.}
  \label{Convergence}
  \label{fig:env}
\end{figure}

\subsection{Rigidity preservation for a leader-follower formation}
The rigidity preservation is generally not the only objective for  mobile robot teams and the proposed controller should be used in conjunction with other controls  to achieve additional tasks. Here we consider a simple and typical leader-follower scenario that has also been considered as simulation examples for the network connectivity control \cite{zavlanos2007potential}, \cite{yang2010decentralized}. Without loss of generality, we suppose a network  of $n=5$ robots in the plane and  the first one is chosen as a leader with additional dynamics. We simply assume that the leader's dynamics is described by $\dot p_1 = u_1^r+ u_1^l$, where $u_1^l = [0.5, 0.3+0.4\text{cos}(p_{1x})]^T$ is an additional input for the leader and is also independent of the followers' motions.

We consider the energy function of (15) with a lower bound $\epsilon=2$ and design the control together with the  distributed estimation algorithm. All the agents are randomly placed such that  the initial formation is infinitesimally rigid. All the followers run the control (16) to preserve the rigidity of the formation as well as to track the motion of the leader. Other simulation settings are the same as above section. The simulation results are shown in Fig.2 and Fig.3.  As revealed in Fig.2, during the leader-follower motion, the rigidity property is always preserved while the formation also allows a flexible geometric shape (the communication edges are allowed to be changed as long as a rigid formation is preserved). From Fig.3 one can observe that the network connectivity and collision avoidance between any two robot  have also been achieved during robots' motions.
\begin{figure}
  \centering
  \includegraphics[width=80mm]{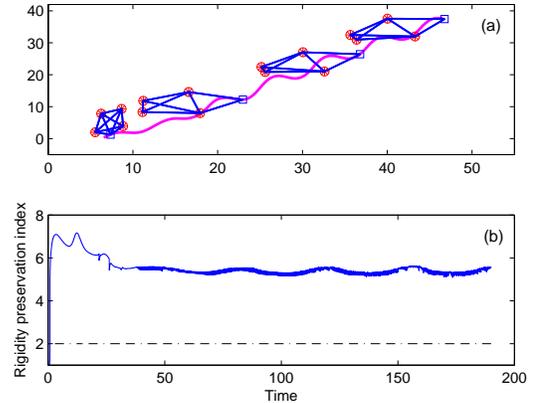}
  \caption{Leader-follower formation with rigidity preservation. (a) Snapshots of the motion. The leader is marked with a square shape and the red line represents the trajectory of the leader; (b) Evolution of $\lambda_4(E)$. The black dashed line represents the lower bound of the rigidity index. }
  \label{fig:env}
\end{figure}

\begin{figure}
  \centering
  \includegraphics[width=80mm]{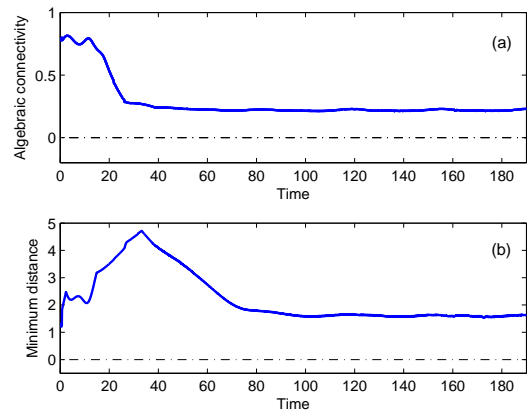}
  \caption{Leader-follower formation with rigidity preservation (cont.). (a) Evolution of the algebraic connectivity eigenvalue; (b) The minimum inter-agent distance during the rigidity preservation and leader tracking process.   }
  \label{Convergence}
  \label{fig:env}
\end{figure}

\section{Concluding remarks}
In this paper we have discussed a rigidity preservation problem by defining and analyzing a rigidity Laplacian matrix. As a natural extension of the standard graph Laplacian, this matrix displays several interesting and useful properties. The rigidity preservation problem is formulated as an eigenvalue control problem, and a gradient control scheme is derived from the defined potential function. To implement the controller in a distributed way, we devise an eigenvector and eigenvalue estimation control via the inverse power iteration method. The results and effectiveness of the distributed algorithms are validated by simulations. Directions for future work include generalizations on more complex robot models, controller design for other coordination objectives (in conjunction with the rigidity preservation scheme in this paper), and experiments on real robots for testing the  algorithms.
\section*{ACKNOWLEDGMENT}

\begin{small}
This work was supported by NICTA, which is funded by the Australian Government
as represented by the Department of Broadband, Communications and the Digital
Economy and the Australian Research Council (ARC) through the ICT Centre of Excellence
program. B. D. O. Anderson was also supported by the ARC under grant DP110100538.  C. Yu was supported by
the Australian Research Council through a Queen Elizabeth II Fellowship and
Discovery Projects DP-110100538 and DP-130103610, and the Overseas Expert
Program of Shandong Province. Z. Sun was also supported by the Prime Minister's Australia Asia Incoming Endeavour Postgraduate Award.
\end{small}

\addtolength{\textheight}{-12cm}   

\bibliographystyle{ifacconf.bst}
\bibliography{rigidity}

\begin{thebibliography}{10}

\bibitem{stump2008connectivity}
E.~Stump, A.~Jadbabaie, and V.~Kumar,
\newblock ``Connectivity management in mobile robot teams'',
\newblock in {\em Robotics and Automation (ICRA) 2008. IEEE International
  Conference on}. IEEE, 2008, pp. 1525--1530.

\bibitem{zavlanos2011graph}
M.~M. Zavlanos, M.~Egerstedt, and G.~J. Pappas,
\newblock ``Graph-theoretic connectivity control of mobile robot networks'',
\newblock {\em Proceedings of the IEEE}, vol. 99, no. 9, pp. 1525--1540, 2011.

\bibitem{sabattin2013}
L.~Sabattini, C.~Secchi, N.~Chopra, and A.~Gasparri,
\newblock ``Distributed control of multirobot systems with global connectivity
  maintenance'',
\newblock {\em Robotics, IEEE Transactions on}, vol. PP, no. 99, pp. 1--6,
  2013.

\bibitem{anderson2008rigid}
B.~D.~O. Anderson, C.~Yu, B.~Fidan, and J.~Hendrickx,
\newblock ``Rigid graph control architectures for autonomous formations'',
\newblock {\em Control Systems, IEEE}, vol. 28, no. 6, pp. 48--63, 2008.

\bibitem{aspnes2006theory}
J.~Aspnes, T.~Eren, D.~K. Goldenberg, A.~S. Morse, W.~Whiteley, Y.~R. Yang,
  B.~D.~O. Anderson, and P.~N. Belhumeur,
\newblock ``A theory of network localization'',
\newblock {\em Mobile Computing, IEEE Transactions on}, vol. 5, no. 12, pp.
  1663--1678, 2006.

\bibitem{kim2010optimizing}
Y.~Kim, G.~Zhu, and J.~Hu,
\newblock ``Optimizing formation rigidity under connectivity constraints'',
\newblock in {\em Decision and Control (CDC), 2010 49th IEEE Conference on}.
  IEEE, 2010, pp. 6590--6595.

\bibitem{zhu2011link}
G.~Zhu and J.~Hu,
\newblock ``Link resource allocation for maximizing the rigidity of multi-agent
  formations'',
\newblock in {\em Decision and Control and European Control Conference
  (CDC-ECC), 2011 50th IEEE Conference on}. IEEE, 2011, pp. 2920--2925.

\bibitem{shames2009minimization}
I.~Shames, B.~Fidan, and B.~D.~O. Anderson,
\newblock ``Minimization of the effect of noisy measurements on localization of
  multi-agent autonomous formations'',
\newblock {\em Automatica}, vol. 45, no. 4, pp. 1058--1065, 2009.

\bibitem{zelazo2013decentralized}
D.~Zelazo, A.~Franchi, H.~H. B{\"u}lthoff, and P.~R. Giordano,
\newblock ``Decentralized rigidity maintenance control with range-only
  measurements for multi-robot systems'',
\newblock {\em arXiv preprint arXiv:1309.0535}, Sep. 2, 2013.

\bibitem{yang2010decentralized}
P.~Yang, R.~A. Freeman, G.~J. Gordon, K.~M. Lynch, S.~S. Srinivasa, and
  R.~Sukthankar,
\newblock ``Decentralized estimation and control of graph connectivity for
  mobile sensor networks'',
\newblock {\em Automatica}, vol. 46, no. 2, pp. 390--396, 2010.

\bibitem{sabattini2012decentralized}
L.~Sabattini, C.~Secchi, and N.~Chopra,
\newblock ``Decentralized connectivity maintenance for networked lagrangian
  dynamical systems'',
\newblock in {\em Robotics and Automation (ICRA) 2012, IEEE International
  Conference on}. IEEE, 2012, pp. 2433--2438.

\bibitem{mesbahi2010graph}
M.~Mesbahi and M.~Egerstedt,
\newblock {\em Graph theoretic methods in multiagent networks},
\newblock Princeton University Press, 2010.

\bibitem{kim2006maximizing}
Y.~Kim and M.~Mesbahi,
\newblock ``On maximizing the second smallest eigenvalue of a state-dependent
  graph laplacian'',
\newblock {\em Automatic Control, IEEE Transactions on}, vol. 51, no. 1, pp.
  116--120, 2006.

\bibitem{zavlanos2007potential}
M.~M. Zavlanos and G.~J. Pappas,
\newblock ``Potential fields for maintaining connectivity of mobile networks'',
\newblock {\em Robotics, IEEE Transactions on}, vol. 23, no. 4, pp. 812--816,
  2007.

\bibitem{ji2007distributed}
M.~Ji and M.~Egerstedt,
\newblock ``Distributed coordination control of multiagent systems while
  preserving connectedness'',
\newblock {\em Robotics, IEEE Transactions on}, vol. 23, no. 4, pp. 693--703,
  2007.

\bibitem{sabattini2011distributed}
L.~Sabattini, N.~Chopra, and C.~Secchi,
\newblock ``Distributed control of multi-robot systems with global connectivity
  maintenance'',
\newblock in {\em Intelligent Robots and Systems (IROS), 2011 IEEE/RSJ
  International Conference on}. IEEE, 2011, pp. 2321--2326.

\bibitem{parlett1980symmetric}
B.~N. Parlett,
\newblock {\em The symmetric eigenvalue problem}, vol.~7,
\newblock SIAM, 1980.

\bibitem{williams2013locally}
R.~K. Williams and G.~S. Sukhatme,
\newblock ``Locally constrained connectivity control in mobile robot
  networks'',
\newblock in {\em Robotics and Automation (ICRA) 2013, IEEE International
  Conference on}, 2013.

\bibitem{freeman2006stability}
R.~A. Freeman, P.~Yang, and K.~M. Lynch,
\newblock ``Stability and convergence properties of dynamic average consensus
  estimators'',
\newblock in {\em Decision and Control, 2006 45th IEEE Conference on}. IEEE,
  2006, pp. 338--343.

\bibitem{bertsekas1997parallel}
D.~P. Bertsekas and J.~N. Tsitsiklis,
\newblock {\em Parallel and Distributed Computation: Numerical Methods},
\newblock Athena Scientific, 1997.

\end{thebibliography}

\end{document}